\documentclass[12pt]{amsart}

\usepackage{amsmath,amsthm,amssymb, mathrsfs,upgreek}
\usepackage{longtable}

\usepackage[margin=2.1cm]{geometry}
\newenvironment{acknowledgments}{\bigskip\small\noindent\emph{Acknowledgments.}}{\par}

\usepackage{graphicx}

\makeatletter
\newcommand*{\da@rightarrow}{\mathchar"0\hexnumber@\symAMSa 4B }
\newcommand*{\da@leftarrow}{\mathchar"0\hexnumber@\symAMSa 4C }
\newcommand*{\xdashrightarrow}[2][]{%
\mathrel{%
\mathpalette{\da@xarrow{#1}{#2}{}\da@rightarrow{\,}{}}{}%
}%
}
\newcommand{\xdashleftarrow}[2][]{%
\mathrel{%
\mathpalette{\da@xarrow{#1}{#2}\da@leftarrow{}{}{\,}}{}%
}%
}
\newcommand*{\da@xarrow}[7]{%
\sbox0{$\ifx#7\scriptstyle\scriptscriptstyle\else\scriptstyle\fi#5#1#6\m@th$}%
\sbox2{$\ifx#7\scriptstyle\scriptscriptstyle\else\scriptstyle\fi#5#2#6\m@th$}%
\sbox4{$#7\dabar@\m@th$}%
\dimen@=\wd0 %
\ifdim\wd2 >\dimen@
\dimen@=\wd2 %
\fi
\count@=2 %
\def\da@bars{\dabar@\dabar@}%
\@whiledim\count@\wd4<\dimen@\do{%
\advance\count@\@ne
\expandafter\def\expandafter\da@bars\expandafter{%
\da@bars
\dabar@ 
}%
}%
\mathrel{#3}%
\mathrel{%
\mathop{\da@bars}\limits
\ifx\\#1\\%
\else
_{\copy0}%
\fi
\ifx\\#2\\%
\else
^{\copy2}%
\fi
}%
\mathrel{#4}%
}
\makeatother

\def\act{\curvearrowright}

\usepackage{extarrows } 
\usepackage{boldline}
\usepackage{mathtools}

\newcommand{\comp}{\mathbin{\scriptstyle{\circ}}}

\newcommand \MMM {\mathscr{M}}

\newcommand{\GQ}{G\mathbb{Q}}
\newcommand{\KK}{\mathbb{K}}
\newcommand \CC {\mathbb{C}}
\newcommand \ZZ {\mathbb{Z}}
\newcommand \PP {\mathbb{P}}

\newcommand \QQ {\mathbb{Q}}
\newcommand{\Aff}{\mathbb{A}}

\newcommand{\ed}{{\operatorname{ed}}}
\newcommand{\PSL}{{\operatorname{PSL}}}

\newcommand{\PGL}{{\operatorname{PGL}}}
\newcommand{\PSp}{{\operatorname{PSp}}}

\newcommand{\GL}{{\operatorname{GL}}}

\newcommand{\jj}{{\operatorname{j}}}
\newcommand{\jjw}{{\bar{\operatorname{j}}}}
\newcommand{\bb}{{\operatorname{b}}}
\newcommand{\g}{{\operatorname{g}}}
\newcommand{\mult}{{\operatorname{mult}}}
\newcommand{\bir}{\underset{\mathrm{bir}}{\sim}}
\newcommand{\Cr}{\operatorname{Cr}}

\newcommand{\Bir}{\operatorname{Bir}}
\newcommand{\Aut}{\operatorname{Aut}}
\newcommand{\mumu}{{\boldsymbol{\mu}}}
\newcommand{\rk} {\operatorname{rk}}
\newcommand{\Pic}{\operatorname{Pic}}
\newcommand{\bF}{\mathbf{F}}
\newcommand{\Sym}{\mathfrak{S}}
\newcommand{\Alt}{\mathfrak{A}}
\newcommand{\Fix}{{\operatorname{Fix}}}
\newcommand{\Fixnu}{{\operatorname{F}^\mathrm{nu}}}

\newcounter{NN}\numberwithin{NN}{subsection}
\renewcommand{\theNN}{\arabic{NN}${}^o$}
\def\nr{\refstepcounter{NN}{\theNN}}

\newcounter{acout}[subsection]

\newtheorem{theorem}[subsection]{Theorem}
\newtheorem{theorem1}[acout]{Theorem}

\newtheorem*{corollary*}{Corollary}

\newtheorem{proposition}[subsection]{Proposition}

\newtheorem{conjecture}[subsection]{Conjecture}

\theoremstyle{definition}

\newtheorem*{definition*}{Definition}

\newtheorem*{example*}{Example}
\newtheorem*{examples*}{Examples}
\newtheorem{problem}[subsection]{Problem}
\newtheorem{problem1}[acout]{Problem}

\newtheorem*{notation*}{Notation}

\numberwithin{equation}{section}

\usepackage[hyperfootnotes=false,colorlinks=true,allcolors=blue]{hyperref}

\begin{document}

\title{{Finite groups of birational transformations}}

\author{Yuri Prokhorov}

\date{}

\address{Steklov Mathematical Institute,
8 Gubkina street, Moscow 119991, Russia \&
\newline
AG Laboratory, HSE, 6 Usacheva str., Moscow, Russia; }
\email{prokhoro@mi-ras.ru }

\begin{abstract}
We survey new results on finite groups of birational transformations of algebraic varieties. 

\end{abstract}
\maketitle

\renewcommand{\thefootnote}{}
\footnotetext{\bf Proceedings of the 8th European Congress of Mathematics}

\section{Introduction}
We work over a field $\Bbbk$ of characteristic $0$.
Typically, unless otherwise mentioned, we assume that 
$\Bbbk$ is algebraically closed. 
The \emph{Cremona group} $\Cr_{n}(\Bbbk)$ of rank $n$ 
is the group of $\Bbbk$-automorphisms of the field $\Bbbk(x_1,\dots,x_n)$ of 
rational functions in $n$ independent variables. Equivalently, 
$\Cr_{n}(\Bbbk)$ can be viewed as the group of birational transformations 
of the projective space $\PP^n$. It is easy to show that for $n=1$ the group 
$\Cr_n(\Bbbk)$ consists of linear projective transformations: 
\[
\Cr_1(\Bbbk)= 
\PGL_2(\Bbbk). 
\]
On the other hand, for $n\ge 2$ the group 
$\Cr_n(\Bbbk)$ has extremely complicated structure. In particular, it contains linear algebraic subgoups 
of arbitrary dimension and has a lot of normal non-algebraic subgroups \cite{CantatLamy, BlancLZ:quotients}.
We refer to \cite{Alberich-Carraminana,Cantat2013,Cantat2018,Deserti:book-s,Hudson1927, Serre-2008-2009} 
for surveys, historical r\'esum\'es,
and introductions to the subject. 

\begin{examples*}
\begin{enumerate}
\item 
Any matrix $A=\|a_{i,j}\|\in \GL_n(\ZZ)$ defines an element $\varphi_A\in \Cr_n(\Bbbk)$ via the following action on $\Bbbk(x_1,\dots,x_n)$:
\[
\varphi_A: x_i \longmapsto x_1^{a_{1,i}} x_2^{a_{2,i}}\cdots x_n^{a_{n,i}}.
\]
Such Cremona transformations are called \emph{monomial}.
For $n=2$ and $A=-\mathrm{id}$ the transformation $\varphi_A$ is known as the \emph{standard quadratic involution}
\[
(x_1,\, x_2)\longmapsto \big(x_1^{-1},\, x_2^{-1}\big).
\]
\item 
Let $S$ be an algebraic variety admitting a generically finite rational map
\[
\pi: S \dashrightarrow \PP^{n-1}
\]
of degree $2$. In an affine piece and suitable coordinates 
$S$ can be given by the equation $y^2=f(x_1,\dots,x_{n-1})$.
One can associate with $(S,\pi)$ an involution $\tau \in \Cr_n(\Bbbk)$ acting on $\Bbbk(x_1,\dots,x_{n-1},y)$ via
\[
\tau: (x_1,\dots,x_{n-1},y) \longmapsto \big(x_1,\dots,x_{n-1},f(x_1,\dots,x_{n-1})\cdot y^{-1}\big). 
\]
If $n=2$ and $S$ is a hyperelliptic curve, then $\tau$ is known as the
\emph{de Jonqui\`eres involution}.
\end{enumerate}
\end{examples*}

The study of the Cremona group has very long history.
Basically, it was started in earlier works of A.~Cayley and L.~Cremona,
and since then this group has
been the object of many studies.
In these notes we concentrate on the following particular problem.

\begin{problem}
Describe the structure of finite subgroups of $\Cr_n(\Bbbk)$.
\end{problem}

Note however that the projective space is not an exceptional variety from 
the algebro-geometric point of view. So one can ask similar question replacing 
$\Cr_n(\Bbbk)$ with the group of birational transformations $\Bir(X)$ of an
arbitrary algebraic variety $X$.
Hence it is natural to pose the following

\begin{problem}
Describe the structure of finite subgroups of 
$\Bir(X)$, where $X$ is an algebraic variety.
\end{problem}
We deal with the most recent results related to these problems. Definitely our survey is not exhaustive.

\section{Equivariant Minimal Model Program}
In this section we collect basic facts on the so-called $G$-Minimal Model 
Program (abbreviated as $G$-MMP). 
This program is the main tool in the study of finite groups of birational 
transformations.
For a detailed exposition we refer to \cite{P:G-MMP}.

Let $G$ be a \emph{finite} group.
Following Yu.~Manin \cite{Manin:raf-surf2e} we say that an algebraic variety $X$ is a 
\emph{$G$-variety} if it is equipped with a regular faithful action $G \act 
X$, i.e. if there exists an injective homomorphism $\alpha: G \hookrightarrow \Aut(X)$.
A morphism (resp. rational map) $f: X\to Y$ of $G$-varieties is a \emph{$G$-morphism} (resp. \emph{$G$-map})
if there exists a group automorphism $\varphi: G\to G$ such that, for any $g\in G$, 
\[
f\comp \alpha(g) = \beta(\varphi(g))\comp f,
\]
where $\alpha: G \hookrightarrow\Aut(X)$ and $\beta: G \hookrightarrow\Aut(Y)$ are the embeddings corresponding to the actions $G \act 
X$ and $G \act 
Y$, respectively.

For any $G$-variety $X$ the action $G \act X$ induces an embedding $G 
\hookrightarrow \Aut_\Bbbk(\Bbbk(X))$ to the automorphism group of the field of 
rational functions $\Bbbk(X) $. Conversely, given any finitely generated extension $\KK/\Bbbk$ and any finite subgroup $G\subset 
\Aut_{\Bbbk}(\KK)$, there exists a $G$-variety $X$ 
and an isomorphism $\Bbbk(X)\simeq_{\Bbbk} \KK $ inducing $G\subset \Aut_{\Bbbk}(\KK)$. Thus, we have.

\begin{proposition}
Let $\KK/\Bbbk$ be finitely generated field extension. Then there exists a 
{\rm 1-1} correspondence 
between finite subgroups $G \subset \Aut_\Bbbk(\KK)$ considered 
modulo conjugacy and $G$-varieties $X$ 
such that
$\Bbbk(X)\simeq_{\Bbbk} \KK $
considered
modulo $G$-birational equivalence.
\end{proposition}

Recall that a variety $X$ is said to be \emph{rational} if it is birationally
equivalent to the projective space $\PP^n$ or, equivalently, if the field extension $\Bbbk(X)/\Bbbk$ is purely transcendental.

\begin{corollary*}
There exists a 
{\rm 1-1} correspondence 
between finite subgroups $G \subset\Cr_n(\Bbbk)$ considered 
modulo conjugacy and rational $G$-varieties $X$ 
such that
$\Bbbk(X)\simeq_{\Bbbk} \KK $
considered
modulo $G$-birational equivalence.
\end{corollary*}

Next, due to equivariant resolution theorem (see e.g.~\cite{Abramovich-Wang})
it is possible to replace $X$ with a smooth projective model.

\begin{proposition}[see e.g. {\cite[14.1.1]{P:G-MMP}}]
\label{prop:reg}
For any $G$-variety $X$ there exists a \emph{smooth projective} $G$-variety 
$Y$ that is $G$-birationally equivalent to $X$. 
\end{proposition} 

Thus the above considerations allow us to reduce the problem 
of classification of finite subgroups of $\Bir(X)$ to 
the study of subgroups in $\Aut(Y)$, where $Y$ is a smooth projective
variety. The main difficulty arising here is that this $G$-variety $Y$ is not unique in its $G$-birational equivalence class.
So, given $G$-birational equivalence class of algebraic $G$-varieties, we need to choose some good representative in it.
This can be done by means of the $G$-MMP.
The higher-dimensional MMP forces us to consider 
varieties with certain very mild, so-called terminal singularities.

\begin{definition*}
A normal variety $X$ has \emph{terminal singularities} if 
some multiple $mK_X$ of the canonical Weil divisor $K_X$ is Cartier and
for any birational morphism
$f: Y\to X$ one can write
\[
mK_Y =f^*mK_X+\sum a_i E_i, 
\]
where $E_i$ are all the exceptional divisors and $a_i>0$ for all $i$. The smallest positive $m$ such that $mK_X$ is Cartier is called the \emph{Gorenstein index} of $X$.
\end{definition*}

\begin{definition*}
A $G$-variety $X$ has \emph{$\GQ$-factorial singularities} if a multiple of any 
$G$-invariant Weil divisor on $X$ is Cartier.
\end{definition*}

It is important to note that terminal 
singularities lie in codimension $\ge 
3$. In particular, terminal surface singularities are smooth.

\begin{example*}[\cite{Mori:term-sing, Reid:YPG}]
Let the cyclic group $\mumu_r$ act on $\Aff^4$ diagonally via
\[
(x_1,x_2,x_3, x_4) \longmapsto \big(\upzeta x_1,\, \upzeta^{-1} x_2,\, \upzeta^a 
x_3, x_4\big),\quad \upzeta=\upzeta_r=\exp (2\uppi \operatorname{i}/r), \quad \gcd(a,r)=1.
\]
Then for a polynomial $f(u,v)$ the singularity of the quotient 
\[
\{x_1x_2+f(x_3^r,x_4)=0\}/\mumu_r
\]
at~$0$ is
terminal whenever it is isolated.
\end{example*}

The aim of the $G$-MMP is to replace a $G$-variety with another one which is ``minimal'' in some sense. As we mentioned above, 
running the $G$-MMP we have to consider singular varieties and the class of terminal $\GQ$-factorial singularities is the smallest class that 
is closed under the $G$-MMP. 

\begin{definition*}[for simplicity we assume that $\Bbbk$ is uncountable] A variety 
$X$ is \emph{uniruled} if
for a general point $x\in X$ there exists 
a rational curve $C\subset X$ passing through $x$.
A variety 
$X$ is \emph{rationally connected} if 
two general points $x_1,\, x_2\in X$ can be connected by a rational curve.
\end{definition*}

Note that a rationally connected surface is rational, and an
uniruled surface is birationally equivalent to $C\times \PP^1$, where $C$ is a curve. 

\begin{definition*}
Let $Y$ be a $G$-variety with only terminal $\GQ$-factorial singularities
and let $f: Y\to Z$ be a $G$-equivariant morphism with connected fibers to a lower-dimensional variety $Z$, where the action of $G$ on $Z$ is not necessarily faithful. Then $f$ is called \emph{$G$-Mori fiber space} (abbreviated as $G$-Mfs) if the anti-canonical class $-K_Y$ is 
$f$-ample
and 
$\rk \Pic(Y/Z)^G=1$. If $Z$ is a point, then $-K_Y$ is 
ample and $Y$ is called \emph{$\GQ$-Fano variety}. 
Two-dimensional $\GQ$-Fano varieties are traditionally called \emph{$G$-del Pezzo surfaces}.
\end{definition*}

\begin{definition*}
A $G$-variety 
$Y$ is said to be a
\emph{$G$-minimal model} if it has only terminal $\GQ$-factorial 
singularities and the canonical class $K_Y$ is numerically effective (nef).
\end{definition*}
It is not difficult to show that the concepts of $G$-minimal model and
$G$-Mori fiber space are mutually exclusive. 
Moreover, if $f: Y\to Z$ is a $G$-Mfs, then its general fiber is rationally connected, hence $Y$ is uniruled. 
On the other hand, a $G$-minimal model is never uniruled \cite{MiyaokaMori}.
The following assertions
are usually formulated for varieties without group actions.
The corresponding equivariant versions can be easily deduced from non-equivariant
ones (see \cite{P:G-MMP}).

\begin{theorem}[{\cite{BCHM}}]
Let $X$ be an uniruled $G$-variety. Then there exists a birational $G$-map
$X \dashrightarrow Y$, where $Y$ has a structure of $G$-Mfs\quad $f: Y\to Z$.
\end{theorem}

\begin{conjecture}
Let $X$ be a non-uniruled $G$-variety. Then there exists a birational $G$-map
$X \dashrightarrow Y$, where $Y$ is a $G$-minimal model.
\end{conjecture}
The conjecture is known to be true in dimension $\le 4$ \cite{Mori:flip, Shokurov:PL:e}, as well as in the case where $K_X$ is big \cite{BCHM}, and in some other cases. In arbitrary dimension a weaker notion of quasi-minimal models works quite satisfactory \cite{Prokhorov-Shramov-J}.

\section{Cremona group of rank $2$}
The $G$-MMP for surfaces is much more simple than in higher dimensions.
It was developed in works of Yu.~Manin and V.~Iskovskikh (see \cite{Manin:raf-surf2e}).
In the two-dimensional case the $G$-MMP works in the category of smooth $G$-surfaces and all the birational transformations are contractions of disjoint unions of $(-1)$-curves.
For a $G$-Mfs
\mbox{$f: Y\to Z$} there are two possibilities:
\begin{enumerate} 
\item \label{2dP}
$Z$ is a point and then $Y$ is a $G$-del Pezzo surface,

\item \label{2cb}
$Z$ is a curve, any fiber of $f$ is a reduced plane conic and $\rk \Pic(Y)^G=2$. In this case $f$ is called \emph{$G$-conic bundle}. 
\end{enumerate}

Thus to study finite subgroups of $\Cr_2(\Bbbk)$ one has to consider 
the above two classes of $G$-Mfs's in detail.
The classification of del Pezzo surfaces is well known and very short. 
Hence, to study the case~\ref{2dP} one has to list all finite subgroups $G\subset \Aut(Y)$ 
satisfying the condition $\rk \Pic(Y)^G=1$.
The full list was obtained Dolgachev and Iskovskikh \cite{Dolgachev-Iskovskikh}.
In contrast, the class of conic bundles is huge and consists of an infinite number of families. In this case a reasonable approach is to find an algorithm
of enumerating conic bundles $Y/Z$ together with subgroups 
$G\subset \Aut(Y/Z)$ satisfying $\rk \Pic(Y)^G=2$.
This also was done by Dolgachev and Iskovskikh \cite{Dolgachev-Iskovskikh}
(see also \cite{Tsygankov:MSb:e}). However even using this algorithm it is very hard to get a complete list of corresponding groups.

As an example, we present well-known classical result on the classification of 
subgroups of order $2$ in $\Cr_2(\Bbbk)$.
It was obtained by E.~Bertini \cite{Bertini1877} in 1877, however his arguments were incomplete from modern point of view. 
A new rigorous proof was given by L.~Bayle and A.~Beauville \cite{Bayle-Beauville-2000}.

\begin{theorem}
\label{thm:invol}
Let $G=\{1,\, \tau\} \subset \Cr_2(\Bbbk)$ be a subgroup of order $2$.
Then the embedding $G\subset \Cr_2(\Bbbk)$ is induced by one of the following 
actions on a rational surface $X$
\par\medskip\noindent
\renewcommand\arraystretch{1.5}{\rm
\begin{longtable}{l|p{0.21\textwidth}|p{0.64\textwidth}}
&\multicolumn{1}{c|}{$\tau$} & \multicolumn{1}{c}{$X$ and $\tau$} 
\\\hlineB{3}
\endhead
&\multicolumn{1}{c|}{$\tau$} & \multicolumn{1}{c}{$X$ and $\tau$} 
\\\hlineB{3}
\endfirsthead
\nr\label{linear}&linear involution & $\PP^2$
\\\hline
\nr\label{de-Jonquieres}&de Jonqui\`eres involution of genus $g\ge 1$& 
$X=\{y_1y_2=p(x_1,x_2)\}\subset \PP(1,1,g+1,g+1)$
\newline
$p$ is a homogeneous form of degree $2g+2$,
\newline
$\tau$ is the deck involution of the projection
\newline $X\xdashrightarrow{\ 2:1\ }\PP(1,1,g+1)$, 
\newline 
$(x_1,x_2,y_1, y_2)\mapsto (x_1,x_2,y_1+y_2)$
\\\hline
\nr\label{Geiser}&Geiser involution&
$X=\{y^2= p(x_1,x_2,x_3)\}\subset \PP(1,1,1,2)$,
\newline
$p$ is a homogeneous form of degree $4$,
\newline
$\tau$ is the deck involution of the projection\newline $X\xlongrightarrow{\ 2:1\ } \PP(1,1,1)=\PP^2$
\\\hline
\nr\label{Bertini}&Bertini involution& 
$X=\{z^2=p(x_1,x_2,y)\}\subset \PP(1,1,2,3)$,
\newline
$p$ is a quasihomogeneous form of degree $6$,
\newline 
$\tau$ is the deck involution of the projection\newline $X\xlongrightarrow{\ 2:1\ } \PP(1,1,2)$
\end{longtable}
}
\par\medskip\noindent
Here $\PP(w_1,\dots,w_n)$ denotes the weighted projective space 
with corresponding weights.
\end{theorem}
In the cases~\ref{linear}, \ref{Geiser}, and~\ref{Bertini} the variety $X$ is a del Pezzo surface of degree $9$, $2$, and $1$, respectively.
In the case~\ref{de-Jonquieres} the projection $X\dashrightarrow \PP(1,1)=\PP^1$ becomes a $G$-conic bundle after blowing up the indeterminacy points. 

The $G$-MMP was successfully applied for classification of various classes 
of finite subgroups in $\Cr_2(\Bbbk)$: groups of prime order \cite{Fernex:Crem}, $p$-elementary groups \cite{Beauville2007}, abelian groups \cite{Blanc:Lin, Blanc:cycl}, and finally arbitrary groups \cite{Dolgachev-Iskovskikh}.
Here is another example of classification results.

\begin{theorem}[{\cite{Dolgachev-Iskovskikh}}]
\label{thm:simple2}
Let $G\subset \Cr_2(\CC)$ be a finite simple group. 
Then $G$ is isomorphic to one of the following:
\begin{equation*}
\Alt_5,\quad \Alt_6,\quad\PSL_2(\bF_7), 
\end{equation*}
where $\Alt_n$ is the alternating group of degree $n$ and $\PSL_n(\bF_q)$
is the projective special linear group over the finite field $\bF_q$. 

Moreover, if $G\not \simeq \Alt_5$, then the embedding $G\subset \Cr_2(\Bbbk)$ is induced by one of the following 
actions on a del Pezzo surface $X$:
\par\smallskip\renewcommand\arraystretch{1.1}
\noindent
\begin{tabular}{l|l|p{0.72\textwidth}}
\multicolumn{1}{c|}{$G$} & \multicolumn{1}{c|}{$|G|$} &\multicolumn{1}{c}{$X$} 
\\\hlineB{3}
$\Alt_6$\quad & $360$\quad& $\PP^2$
\\
$\PSL_2(\bF_7)$& $168$& $\PP^2$
\\
$\PSL_2(\bF_7)$ & $168$& 
$\{y^2=x_1^3x_2+x_2^3x_3+ x_3^3x_1\}\subset \PP(1,1,1,2)$
\end{tabular}
\end{theorem} 
A complete classification of embeddings $\Alt_5\hookrightarrow\Cr_2(\Bbbk)$ can be found in 
\cite{Cheltsov:2ineq}.

\section{Cremona group of rank $3$}

The MMP in dimension $3$ is more complicated than two-dimensional one but still it is developed very well. In particular, terminal threefold singularities 
are classified up to analytic equivalence \cite{Mori:term-sing, Reid:YPG}. The structure of all intermediate steps of the MMP and Mfs’s 
is also studied relatively well (see \cite{P:G-MMP} for a survey). 

For a three-dimensional $G$-Mori fiber space $f: Y\to Z$ there are three possibilities:
\begin{enumerate} 
\item \label{3F}
$Z$ is a point, then $Y$ is a (possibly singular) $\GQ$-Fano threefold;

\item \label{3dPf}
$Z$ is a curve, then $f$ is called a $\GQ$-del Pezzo fibration,

\item \label{3cb}
$Z$ is a surface, then $f$ is a $\GQ$-conic bundle.
\end{enumerate}
A $\GQ$-conic bundle can be birationally transformed to a \emph{standard $G$-conic bundle},
i.e. $\GQ$-conic bundle such that both $X$ and $Z$ are smooth \cite{Avilov:cb}.
For $\GQ$-del Pezzo fibrations there are only some partial results of this type 
(see \cite{Corti:DP, Loginov:dP1}).
Nevertheless, the main difficulty in the application $G$-MMP to the classification of finite groups of birational transformations is the lack of a complete classification of Fano threefolds with terminal 
singularities. 
At the moment only some very particular classes of $\GQ$-Fano threefolds 
are studied (see \cite{IP99,Avilov:cubic:Eu,Avilov-V4e,P:GFano-all,Prokhorov-v22} and references therein).
Some roundabout methods work in the case of 
``large'' in some sense (in particular, simple) finite groups. 

\begin{theorem}[\cite{P:JAG:simple}]
\label{thm:simple}
Let $G\subset \Cr_3(\CC)$ be a finite simple subgroup.
Then $G$ is isomorphic to one of the following:
\begin{equation*}
\Alt_5,\quad \Alt_6,\quad \Alt_7, \quad\PSL_2(\bF_7),\quad \PSL_2(\bF_8),\quad \PSp_4(\bF_3),
\end{equation*}
where $\PSp_4(\bF_3)$ is the projective symplectic group over $\bF_3$.

All the possibilities occur.
\end{theorem}
This classification is a consequence of the following more general result.

\begin{theorem}[{\cite{P:JAG:simple}}]
\label{thm:simple3}
Let $Y$ be a rationally connected threefold and let $G\subset \Bir (Y)$ 
be a finite simple group.
If $G$ is not embeddable to $\Cr_2(\CC)$, then $Y$ is $G$-birationally equivalent to one of the following $\GQ$-Fano threefolds.
\renewcommand\arraystretch{1.5}
\par\medskip\noindent
{\rm\begin{tabular}{l|p{0.11\textwidth}|p{0.61\textwidth}|p{0.10\textwidth}}
& \multicolumn{1}{c|}{$G$} & \multicolumn{1}{c|}{$X$} & rational?
\\\hlineB{3}
\nr\label{tab:A7-V6}
& $\Alt_7$ & 
$X_6'=\{\upsigma_{1, 7}=\upsigma_{2, 7}=\upsigma_{3, 7}=0\}\subset \PP^5\subset \PP^6$& 
NO
\\
\nr\label{tab:A7-P3}
& $\Alt_7$ & $\PP^3$& YES
\\
\nr\label{tab:PSp-P3}
& $\PSp_4(\bF_3)$ & $\PP^3$& YES
\\
\nr\label{tab:PSp-Bur}
& $\PSp_4(\bF_3)$ & Burkhardt quartic
$X_4^{\mathrm b}=\{\upsigma_{1, 6}=\upsigma_{4, 6}=0\}\subset \PP^4\subset \PP^5$&YES
\\
\nr\label{tab:SL2-8}
& $\PSL_2(\bF_8)$ & special Fano threefold $X_{12}^{\mathrm m}\subset \PP^8$
of genus $7$&YES
\\
\nr\label{tab:PSL2-11-Kc}
& $\PSL_2(\bF_{11})$ & Klein cubic $X^{\mathrm k}_3=\{x_1x_2^2+x_2x_3^2+\cdots 
x_5x_1^2=0\}\subset \PP^4$&NO
\\
\nr\label{tab:PSL2-11-X14}
& $\PSL_2(\bF_{11})$ & special Fano threefold $X_{14}^{\mathrm a}\subset \PP^9$ 
of genus $8$&NO
\end{tabular}}
\par\smallskip\noindent
where $\upsigma_{d, k}=\upsigma_{d, k}(x_1,\dots,x_k)$ is the elementary symmetric polynomial of degree $d$ in $k$ variables.
\end{theorem} 

Below we outline the proof of Theorem~\ref{thm:simple3}.

Assume that $G$ is not embeddable to $\Cr_2(\Bbbk)$, i.e. it is not isomorphic to any of the groups listed
in Theorem~\ref{thm:simple2}.
First, Proposition~\ref{prop:reg} allows us to assume that the action of $G$ is regularized
on some smooth projective $G$-variety $X$.
By running the equivariant 
MMP, we may assume that $X$ has a structure of a
$G$-Mfs
$f:X\to Z$ (because $X$ is rationally connected).
Consider the case $\dim Z>0$. Since $G$ is a simple group,
it must act faithfully on the base $Z$ or on the general fiber $F$.
Since the varieties $F$ and $Z$ are rational, this means that $G$ is contained
in the plane Cremona group $\Cr_2 (\Bbbk)$.
The contradiction proves Theorem~\ref{thm:simple3} in
the case $\dim Z>0$.

Hence, we may further assume that
$Z$ is a point and $X$ is a $\GQ$-Fano threefold.
Consider the case where $X$ is not Gorenstein, i.e.
the canonical class $K_X$ is not a Cartier divisor. It turns out that this case
does not occur. 
Let $P_1,\dots,P_n\in X$ be all non-Gorenstein points 
and let $r_1,\dots,r_n$ be the corresponding Gorenstein indices.
Arguments based on Bogomolov-Miyaoka inequality (see \cite{Kawamata:bF, KMMT-2000} and \cite[\S 12]{P:G-MMP}) show that
\[
\textstyle \sum \big(r_i -\frac 1 {r_i}\big) <24.
\]
Hence, $n\le 15$. Then using the classification of 
transitive
actions of simple groups \cite{atlas} and analyzing the action of stabilizers of $P_i$ one obtains the only possibility:
\begin{itemize} 
\item 
$n=11$, $G\simeq \PSL_2(\mathbf F_{11})$, $r_1=\cdots =r_n=2$.
\end{itemize}
This case is excluded by a more detailed
geometric consideration (see \cite[\S~6]{P:JAG:simple}).

Thus, we may assume that 
$K_X$ is a Cartier divisor. 
In this case according to \cite{Namikawa:Fano} the variety $X$ has a smoothing,
that is,
there exists a one-parameter flat family $\mathfrak{X}/\mathfrak{B}\ni o$ such that the special fiber $\mathfrak{X}_o$ is isomorphic to $X$ and a general geometric fiber $\mathfrak{X}_t$ is a smooth Fano threefold.
Hence some discrete invariants of $X$, such as the Picard lattice $\Pic(X)$ and 
the anticanonical degree $-K_X^3$,
are the same as for smooth Fano threefolds which are completely classified (see \cite{IP99}).
Recall that the Fano index $\iota(X)$ of $X$ is the maximal integer that divides the canonical class $K_X$ in the lattice $\Pic(X)$ \cite{IP99}. 
By \cite[Part II]{P:GFano-all}, we have $\rk \Pic(X)\le 4$. Since $\Pic(X)^G\simeq
\ZZ$ and a simple group that is not isomorphic to $\Alt_5$ cannot have a 
nontrivial integer representation of dimension $\le 4$, we have $\rk \Pic (X)=1$. 
If
$\iota(X)\ge 4$ (resp, $\iota(X)=3$), then $X$ is isomorphic to the projective space $\PP^3$ 
(resp. a quadric in
$\PP^4$)
\cite{IP99}. Then 
from
the classification of finite subgroups in $\PSL_4 (\Bbbk)$ and $\PSL_5 (\Bbbk)$ we get cases~\ref{tab:A7-P3} and~\ref{tab:PSp-P3}. 
Three-dimensional Fano varieties with $\iota(X)=2$ are called del Pezzo threefolds. $G$-Fano threefolds of this type were studied in 
\cite[Part I]{P:GFano-all}. As a consequence of these results we get the
case of the group $G=\PSL_2(\bF_{11})$ acting on the Klein cubic
(case~\ref{tab:PSL2-11-Kc}).  

Finally, let $\Pic (X)=\ZZ\cdot K_X$. 
Recall that in this case the anticanonical degree is written in the form $-K_X^3=2\g(X)-2$, where $\g(X)\in \{2,3,\dots,10,12\}$ \cite{IP99}.
For $\g(X)\le 5$ the variety $X$ has a natural
embedding to a (weighted) projective space as a complete intersection
\cite{IP99}. Using this and some facts from representation theory, we obtain
for the group $G$ two cases~\ref{tab:A7-V6} and~\ref{tab:PSp-Bur}.
The case $\g (X)=6$ can be excluded
using \cite[Corollary 3.11]{Debarre-Kuznetsov:GM}. For $\g(X)\ge 7$
the variety $X$ must be smooth (see \cite[Lemma 5.17]{P:JAG:simple} and
\cite{Prokhorov-v22}). 
Further, using some facts about automorphisms 
of smooth
Fano threefolds \cite{KPS:Hilb} we obtain for the group $G$ two possibilities 
\ref{tab:SL2-8} and \ref{tab:PSL2-11-X14}.
This completes our sketch of the proof of Theorem~\ref{thm:simple3}.\qed  

Similar technique was also applied to the study of finite $p$-subgroups 
and quasi-simple subgroups
in $\Cr_3(\Bbbk)$, see 
\cite{P:p-groups, Prokhorov-2-elementary, Prokhorov-Shramov-p-groups, Kuznetsova:Finite3, Loginov:3Cr} and \cite{BCDP:QS}.

Note that Theorem~\ref{thm:simple3} does not 
describe \emph{embeddings} of 
groups $\Alt_5$, $\Alt_6$, and $\PSL_2(\bF_7)$ to the space Cremona group. It is obvious 
that
such embeddings exist, but their full classification should be significantly 
more difficult.
There are only some partial results in this direction (see e.g. 
\cite{Cheltsov-Przyjalkowski-Shramov-Barth-sextic, Cheltsov-Shramov-2012, Cheltsov-Shramov:5, CheltsovShramov:book, Krylov2020FamiliesOE}).  

\section{Jordan property}
The methods and results of \cite{Dolgachev-Iskovskikh} show that 
one cannot expect a reasonable classification of all finite 
subgroups of Cremona groups of higher rank. 
Thus it is natural to concentrate on the study of general properties 
of these subgroups.
Recall the following two famous results by C.~Jordan and H.~Minkowski.

\begin{theorem}[\cite{Jordan:1878}] 
There exists a function $\jj(n)$ such that for any finite subgroup $G \subset 
\GL_n (\CC)$ 
there exists a normal
abelian subgroup $A \subset G$ of index at most $\jj(n)$.
\end{theorem}

\begin{theorem}[\cite{Minkowski:87}]
There exists a function $\bb(n)$ such that for every finite subgroup $G \subset 
\GL_n (\QQ)$ one has $|G|\le \bb(n)$. 
\end{theorem}

J.-P. Serre \cite{Serre:MMJ, Serre:problems} asked if these properties hold for Cremona groups.
Complete answers to these questions were given in \cite{Prokhorov-Shramov-J, ProkhorovShramov-RC} (see below).
The following very convenient definitions were suggested by V.~L.~Popov \cite{Popov:Russel}. 
\begin{definition*}
\begin{itemize}
\item 
A group $\Gamma$ is \emph{Jordan} if there exists a constant $\jj(\Gamma)$ such 
that any finite
subgroup $G \subset \Gamma$ has a normal abelian subgroup $A$ of index $[G : 
A]\le \jj(\Gamma)$. 
\item 
A group $\Gamma$ is \emph{bounded} \textup(or satisty \emph{bfs} property\textup) if there exists a constant $\bb(\Gamma)$ 
such that for any finite subgroup $G \subset \Gamma$ one has $|G|\le 
\bb(\Gamma)$.
\end{itemize}
\end{definition*} 

\subsubsection*{Rationally connected varieties}

\begin{theorem}[\cite{ProkhorovShramov-RC} \& \cite{Birkar-BAB}]
\label{thm:J-RC}
Let $X$ be a rationally connected variety. Then~$\Bir(X)$ is Jordan.
Moreover, $\Bir(X)$ is \emph{uniformly} Jordan, that is, the constant 
$\jj(\Bir(X))$
depends only on $\dim(X)$.
\end{theorem}

As a consequence we obtain that the group 
$\Cr_n(\Bbbk)$ is Jordan.

Originally Theorem~\ref{thm:J-RC} was proved modulo so-called BAB conjecture (in a weak form)
which is now settled by C.~Birkar:

\begin{theorem}[\cite{Birkar-BAB}]
\label{bab}
Fix $d>0$. The set of all Fano varieties $X$ of dimension at most $d$ with at worst 
terminal singularities form a bounded family, i.e. they are parameterized by a 
scheme of finite type. 
\end{theorem}   

It follows from Theorem \ref{thm:J-RC} that
there is a constant $L = L(n)$ such 
that for any rationally connected variety $X$ of dimension $n$ and for any 
prime $p > L(n)$, every finite $p$-subgroup of $\Bir(X)$ is abelian and 
generated by at most $n$ elements (see \cite{ProkhorovShramov-RC}).
Recently this result was essentially improved 
by Jinsong Xu \cite{Xu:p-groups}: he showed that $L(n)=n+1$.
The proof is based on a result by O.~Haution \cite{Haution:19}.
Thus we have  

\begin{theorem} 
\label{thm:p-groups}
Let $X$ be a rationally connected variety of dimension $n$ and let $G\subset 
\Bir(X)$ be a finite $p$-subgroup. 
If $p > {n+1}$, then $G$ is abelian and is generated by at most $n$ elements.
\end{theorem}

The results of Theorems~\ref{thm:J-RC} and~\ref{thm:p-groups} were applied in the proof of Jordan property of local fundamental groups of log terminal singularities \cite{Braun-Filipazzi-Moraga-Svaldi,Moraga:toroidalization}. 

\subsubsection*{Varieties over non-closed fields}

\begin{theorem}[{\cite{Prokhorov-Shramov-J}} \& \cite{Birkar-BAB}]
\label{thm:b}
Let $X$ be a variety 
over a field $\Bbbk$ of characteristic $0$ which is finitely generated over $\QQ$. Then the group~\mbox{$\Bir(X)$} is bfs.
\end{theorem}
Similar to Theorem~\ref{thm:J-RC}, the proof of this result is based on the BAB conjecture (Theorem~\ref{bab}).

In the case $X=\PP^2$ an explicit bound 
was obtained in \cite{Serre:MMJ} (see also \cite{Dolgachev-Iskovskikh:p}) in terms of cyclotomic invariants of the field $\Bbbk$.
Theorem~\ref{thm:b} can be reformulated in an algebraic form which gives the positive answer to a question of J.-P.~Serre 
\cite{Serre:problems}.

\begin{theorem1}
\label{tmm:fgQa}
Let $\KK$ be a finitely generated field over $\QQ$.
Then the group $\Aut (\KK)$ is bfs.
\end{theorem1} 

\subsubsection*{Jordan constants}

Define the \emph{Jordan constant} of a group $\Gamma$
as the number $\jj(\Gamma)$ that appear in the definition of Jordan property.
The \emph{weak Jordan 
constant} $\jjw(\Gamma)$ 
of $\Gamma$ is the minimal $j$ such that 
for any finite
subgroup $G \subset \Gamma$ there exists an abelian (not necessarily 
normal) subgroup
$A \subset G$ such that $[G: A]\le j$.
Easy group-theoretic arguments show that 
\[
\jjw(\Gamma)\le \jj(\Gamma)\le \jjw(\Gamma)^2.
\]
The exact value of the Jordan constant is known only for Cremona group of rank two: $\jj(\Cr_2(\Bbbk)) = 7200$ (see \cite{Yasinsky-J-const}). On the other hand, weak Jordan 
constants are easier to compute. It was proved in \cite{Prokhorov-Shramov-J-const} that 
\[
\jjw(\Cr_2)=288,\quad\text{} \quad 
\jjw(\Cr_3)=10368. 
\]
Moreover, the inequality $\jjw(\Bir(X))\le 10368$ holds for any rationally connected threefold~$X$.  

\subsubsection*{Jordan property of arbitrary varieties}
It turns out that the group of birational transformations of an algebraic variety is not always Jordan.
The first example was discovered by Yu. Zarhin.
\begin{example*}[\cite{Zarhin10}]
Let $C$ be an elliptic curve and let $X=C\times \PP^1$.
Then the group $\Bir(X)$ is not Jordan.
\end{example*}
On the other hand, the exceptions as above are very rare:
\begin{theorem}[V.~L.~Popov \cite{Popov:Russel}]
\label{thm:Popov-Zarhin}
Let $X$ be an algebraic surface. 
The group
$\Bir(X)$ is not Jordan if and only if 
$X$ is birationally equivalent to $\PP^1\times C$, where $C$ is an elliptic 
curve.
\end{theorem}

The proof of this result given in \cite{Popov:Russel} essentially uses 
a result of I.~Dolgachev which in turn is based on the 
classification of algebraic surfaces. 
Later Theorem~\ref{thm:Popov-Zarhin} was generalized to higher dimensions 
with classification independent proofs.

\begin{theorem}[{\cite{Prokhorov-Shramov-J}}]
\label{thm:3J}
Let $X$ be an algebraic variety.
Then the following assertions hold.
\begin{enumerate}
\item\label{thm:3J-1}\label{thm:3J-2}
If $X$ either is non-uniruled or has irregularity
$q(X)=0$, then $\Bir(X)$ is Jordan.

\item\label{thm:3J-3}
If $X$ is
non-uniruled and
$q(X)=0$, then $\Bir(X)$ is bfs.
\end{enumerate}
\end{theorem}
Similar to Theorems~\ref{thm:b} and \ref{thm:J-RC} the proof of \ref{thm:3J}\ref{thm:3J-2} is
based on boundedness of terminal Fano varieties (Theorem~\ref{bab}).

In dimension three there is the following
much more precise result.

\begin{theorem}[\cite{P-Shramov:3folds}]
Let $X$ be a three-dimensional algebraic variety.
Then $\Bir(X)$ is not Jordan
if and only if either
\begin{enumerate}
\item 
$X$ is birationally equivalent to $ C\times\PP^2$,
where $C$ is an elliptic curve, or 
\item 
$X$ is birationally equivalent to $ S\times\PP^1$, where $S$
is one of the following:
\begin{itemize} 

\item
a surface of Kodaira dimension $\varkappa(S)=1$ such that the Jacobian fibration
of the pluricanonical map $\phi\colon S\to B$ is locally trivial;

\item 
$S$ is either 
an abelian or bielliptic surface \textup(and $\varkappa(S)=0$\textup).
\end{itemize}
\end{enumerate}
\end{theorem}

Below we explain the main
idea of the proof of the necessity.
So we assume that $\Bir(X)$ is not Jordan.
By Theorems~\ref{thm:J-RC} and~\ref{thm:3J} the variety
$X$ is uniruled but it is not rationally connected.
Hence 
there exists a map $X\dashrightarrow Z$ with rationally connected fibers (so-called maximal rationally connected fibration) such that $Z$ is not 
uniruled and $\dim(Z)=1$ or $2$ (see \cite{Kollar-Miyaoka-Mori-1992a}). We have a natural exact sequence
\[
1 \longrightarrow \Bir(X_\eta) \longrightarrow \Bir(X) \longrightarrow \Bir(Z),
\]
where $X_\eta$ is the generic scheme-theoretic fiber.
Since $X_\eta$ is rationally connected and $Z$ is not 
uniruled, the groups $\Bir(X_\eta)$ and $\Bir(Z)$ must be Jordan.
Then group-theoretic arguments show that both groups $\Bir(X_\eta)$ and $\Bir(Z)$ are not 
bfs (see e.g.~\cite[Lemma~2.8]{Prokhorov-Shramov-J}). 
In the case where $Z$ is a curve this implies that $Z$ is elliptic and 
applying the following fact with $\KK=\Bbbk(Z)$ and $S:= X_\eta$ we obtain that $X$ is birationally equivalent to $Z\times\PP^2$.

\begin{proposition}[\cite{P-Shramov:3folds}]
Let $\KK$ be a field containing all roots of $1$ and let $S$ be a surface over $\KK$ such that
$S$ is not
$\KK$-rational, 
$S$ is $\bar \KK$-rational, and $S(\KK)\neq \varnothing$.
Then the group $\Bir(S)$ is bfs.
\end{proposition}
Note that the condition of the existence of a $\KK$-point on $S$ in the above statement is important.
The groups of (birational) automorphisms of geometrically rational surfaces without rational points were studied in the series of papers
\cite{Shramov:aut-2cubics,Shramov:SB-bir,Shramov:SBf}. 

Now assume that $Z$ is a surface. 
According to the main result of \cite{BandmanZarhin:17} the threefold $X $ is birationally equivalent to $Z\times \PP^1$.
By Theorem~\ref{thm:3J} we have $q(Z)>0$.
Thus in the case $\varkappa(Z)=0$ the surface $Z$ must be either abelian or bielliptic.
Since the group $\Bir(Z)$ is not finite in our case, $Z$ cannot be a surface of general type. Consider the case 
$\varkappa(Z)=1$. Then the pluricanonical map $\phi: Z 
\longrightarrow B$ is a $\Bir(Z)$-equivariant elliptic fibration. 
Let 
\[
\mathrm{Jac}(\phi) :E \longrightarrow B
\]
be the corresponding Jacobian fibration.
The automorphism group $\Aut(Z_\eta)$ of the generic fiber $Z_\eta$ over $B$ 
is embedded to $\Bir(Z)$ as a normal subgroup.
Analyzing singular fibers one can conclude that $\Aut(Z_\eta)$ is of finite index in $\Bir(Z)$. In turn, 
$\Aut(Z_\eta)$ has a subgroup $\Aut'(Z_\eta)$ of index at most $6$ isomorphic to the group of $\Bbbk(B)$-points of $E_\eta$. 
Assume that the fibration $\mathrm{Jac}(\phi)$ is not locally trivial.
Then by the functional version of Mordell-Weil theorem, known as Lang-N\'eron theorem; see e.g. \cite{Conrad:Lang-Neron}, 
the group of $\Bbbk(B)$-points of $E_\eta$
is finitely generated, and in particular the torsion subgroup of the group of
points of $E_\eta$ is finite. This implies that $\Aut'(Z_\eta)$ is finite.
\qed

\section{Invariants and rigidity}
The most important part of the classification of finite subgroups in $\Bir(X)$
is to distinguish conjugacy classes.
\begin{problem}
Let $G,\, G'\subset \Bir(X)$ be finite subgroups such that $G\simeq G'$.
How can one conclude that $G$ and $G'$ are \emph{not} conjugate?
\end{problem}
This is equivalent to the following
\begin{problem1}
Let $X$ and $X'$ be $G$-varieties. 
How can one conclude that $X$ and $X'$ are \emph{not} $G$-birational?
\end{problem1}

Below we describe a few approaches to solve the above problems.
Note however that there are no universal methods. 

\subsubsection*{Fixed point locus}
Let $X$ be a smooth projective $G$-variety. 
By $\Fix(X,G)$ we denote the set of $G$-fixed points. 
It is not difficult to show (see \cite{P:invol}) that $\Fix(X,G)$ has at most 
one codimension one component that is not uniruled.
Denote this component by $\Fixnu(X,G)$. 
This is a natural birational invariant in the category of smooth projective $G$-varieties.

\begin{proposition}[\cite{P:invol}]
Let $X$ and $X'$ be smooth projective $G$-varieties. 
If $X$ and $X'$ are $G$-birational, then $\Fixnu(X,G_0)$ and 
$\Fixnu(X',G_0)$ are birational for any subgroup $G_0\subset G$.
\end{proposition}

If $G_0\subset G$ is a normal subgroup, then the set $\Fixnu(X,G_0)$ (if it is not empty) has a structure 
of $(G/G_0)$-variety.
Clearly, the birational type of this $(G/G_0)$-variety is also a birational invariant (cf. \cite{Blanc:cycl}).

\begin{example*} 
According to Theorem~\ref{thm:invol} for subgroups $G\subset \Cr_2(\Bbbk)$ 
of order $2$ we have one of the following possibilities:
\par\medskip\noindent
\renewcommand\arraystretch{1.1}
\begin{center}
\begin{tabular}{l|p{0.37\textwidth}|p{0.49\textwidth}}
&\multicolumn{1}{c|}{involution $\tau \in G$} & 
\multicolumn{1}{c}{$\Fixnu(X,G)$} 
\\\hlineB{3}
\ref{linear}&linear on $\PP^2$ & $\varnothing$
\\
\ref{de-Jonquieres}&de Jonqui\`eres of genus $g\ge 1$& hyperelliptic curve of genus $g$
\\
\ref{Geiser}&Geiser & non-hyperelliptic curve of genus $3$
\\
\ref{Bertini}&Bertini& special non-hyperelliptic curve of genus $4$
\end{tabular}
\end{center}
\end{example*}
Thus the curve $\Fixnu(X,G)$ distinguishes conjugacy classes in this case. The same assertion is true for subgroups of prime order \cite{Fernex:Crem}
but it fails in general \cite{Blanc:Lin}.

\subsubsection*{Cohomological invariants} 
It is not difficult to see that for a smooth 
projective $G$-variety $X$
the cohomology group 
\[
H^1(G,\ \Pic(X)) 
\]
is a $G$-birational invariant (see \cite{Bogomolov-Prokhorov}).
More generally, we say that $G$-varieties $X$ and $X'$ are \emph{stably}
$G$-birationally equivalent if for some $n$ and $m$ the products $X\times 
\PP^n$ and 
$X'\times \PP^m$ are $G$-birationally equivalent, where the action of $G$ on 
$\PP^n$ and $\PP^m$ is supposed to be trivial. Then we have.

\begin{theorem}[\cite{Bogomolov-Prokhorov}]
\label{thm:H1a}
Let $X$ and $X'$ be smooth projective $G$-varieties. If $X$ and 
$X'$ are stably $G$-birationally equivalent, then 
\[
H^1(G,\Pic(X))\simeq H^1(G,\ \Pic(X')).
\]
\end{theorem}
Surprisingly, in some cases the invariant $H^1(G,\Pic(X))$ can be computed in terms of $G$-fixed locus.

\begin{theorem}[\cite{Bogomolov-Prokhorov}]
\label{thm:H1}
Let $G$ be a cyclic group of prime order $p$ and let $X$ be a smooth projective rational $G$-surface. Assume that $\Fixnu(X,G)$ is
a curve of genus $g$. Then
\[
H^1(G, \Pic(X ))\simeq (\ZZ/p \ZZ)^{2g}. 
\]
\end{theorem}
This theorem was slightly generalized with more conceptual proof in \cite{Shinder}.
Another cohomological invariant which is called \emph{Amitsur group} was 
introduced in \cite{BCDP:QS}. 

As a consequence of Theorem~\ref{thm:H1}, one can see that involutions from different families in Theorem~\ref{thm:invol} are not stably conjugate in $\Cr_2(\Bbbk)$.
Note however, that $H^1(G, \Pic(X ))$ is a discrete invariant. 
For example, stable conjugacy of involutions whose 
curves $\Fixnu(X,G)$ are non-isomorphic but have 
the same genus is not known. 

A natural question that arises here is to find examples of subgroups in $\Cr_n(\Bbbk)$ that are stably conjugate but not conjugate.
This question is similar to the birational Zariski problem \cite{BCTSSD85}.

\begin{example*}
Let $G=\Sym_3\times \mumu_2$. 
There are two 
embeddings of this group into the Cremona group $\Cr_2(\Bbbk)$
induced by the following actions:
\begin{enumerate}
\item 
action on $\PP^2=\{x_1 + x_2 + x_3 = 0\}\subset \PP^3$ by permutation and reversing
signs; 
\item 
action on the sextic del Pezzo surface $\{y_1y_2y_3 =y_1'y_2'y_3' \}\subset \PP^1\times\PP^1\times \PP^1$ by permutation and taking inverses.
\end{enumerate}
It was shown in \cite{Lemire-Popov-Reichstein} that these two subgroups in $\Cr_2(\Bbbk)$ are stably conjugate, in fact, they are conjugate in $\Cr_{4}(\Bbbk)$. On the other hand, they are not conjugate \cite{Isk:Two1e}. 
\end{example*}  

Here is another example of this kind wich was pointed out to us by Yuri Tschinkel.

\begin{example*}[{\cite{Reichstein-Youssin:bir}}]
Let $V$ and $W$ be faithful linear representations of $G$ with $\dim(V)=\dim(W)=n$. 
Assume that the images of $G$ in $\GL(V)$ and $\GL(W)$  do not contain non-identity scalar matrices.
Then by a variant of the no-name lemma \cite{Dolgachev:fields-inv} we have the following $G$-birational equivalences of $G$-varieties: 
\[
\PP(V)\times \Bbbk^{n+1} \bir   V\times W \bir \PP(W)\times \Bbbk^{n+1}
\]
where $\Bbbk^{n+1}$ is viewed as the trivial representation. Hence $G$-varieties $V$ and $W$ are stably $G$-birationally equivalent. On the other hand, it may happen that 
they are \emph{not} $G$-birationally equivalent. 

For example, Reichstein and Youssin \cite{Reichstein-Youssin:bir} showed that the {\em determinant} of the action in the tangent space at a fixed point of a finite abelian group, up to sign, is a birational invariant of the action. This allowed them to produce nonbirational linear actions, e.g., of groups $\mumu_{p^n}$ on $\PP^n$, with $p\ge 5$. Many new examples of nonbirational linear actions were given in \cite[Sect.~10-11]{Kresch-Tschinkel:Burnside-repr}; these are based on new invariants introduced in \cite{Kresch-Tschinkel:Burnside-vol} (see also \cite{kontsevich-Pestun-Tschinkel, Hassett-Kresch-Tschinkel:symbols}).
These invariants take into account more refined information about the action on subvarieties with nontrivial abelian stabilizers. 
\end{example*}

A prime number $p$ is said to be a \emph{torsion prime} for the group $\Bir(X)$ if there is a finite abelian $p$-subgroup $G\subset \Bir(X)$ not contained in any algebraic torus of $\Bir(X)$ \cite{Popov-tori}.
Note that if a group $G$ is contained in an algebraic torus 
$T\subset \Bir(X)$, then for any smooth projective birational model $Y$ of $X$ 
on which $T$ acts biregularly we have $H^1(G, \Pic(Y))=0$.
Then by Theorem~\ref{thm:H1a} the inequality 
$H^1(G, \Pic(Y))\neq 0$ for a finite $p$-subgroup $G\subset \Aut(Y)$ implies that a prime number $p$ is a torsion prime for $\Bir(Y)$ and for $\Bir(Y\times \PP^n)$ for any $n$. Using Theorem~\ref{thm:H1} and the classification \cite{Fernex:Crem} one can immediately see
that the set of all torsion primes for $\Cr_2(\Bbbk)$ 
is equal to $\{2,\, 3,\, 5\}$
and the numbers
$2$, $3$, and $5$ 
are torsion primes for $\Cr_n(\Bbbk)$ for any $n\ge 2$. This fact was proved in \cite{Popov-tori} by using another arguments.
In the case $n\ge 3$ the collection of all torsion primes for $\Cr_n(\Bbbk)$
is unknown. 
\subsubsection*{Maximal singularities method}

Maximal singularities method 
is the most powerful tool to study 
birational maps between Mfs's. 
It goes back to works of G.~Fano and even earlier works of other Italian geometers. However the first application of this techniques with rigorous proofs appeared much later 
in the breakthrough paper of Manin and Iskovskikh \cite{Isk-Manin}.
For an introduction to the ``standard'', non-equivariant maximal singularities method we refer to the book \cite{Pukhlikov:book}.
Below we outline very briefly an equivariant version of the method. 

\begin{definition*}[{\cite[Definition~7.10]{Dolgachev-Iskovskikh}},  {\cite[Definition~3.1.1]{CheltsovShramov:book}}]
A $\GQ$-Fano variety $X$ is said to be \emph{$G$-birationally rigid} if 
given birational $G$-map $\Phi:
X\dashrightarrow X^{\sharp}$ to the total space of another $G$-Mfs\ $X^{\sharp}/Z^{\sharp}$, there exists a
birational $G$-selfmap $\psi: X\dashrightarrow X$ such that the composition
$\Phi\circ\psi: X\dashrightarrow X^{\sharp}$ is an isomorphism (in particular, $Z^{\sharp}$ is a point, i.e. $X^{\sharp}$ is also a $\GQ$-Fano variety).

A $\GQ$-Fano variety $X$ is said to be \emph{$G$-birationally superrigid} if any birational $G$-map $\Phi:
X\dashrightarrow X^{\sharp}$ to the total space of another $G$-Mfs\ $X^{\sharp}/Z^{\sharp}$ is an isomorphism.
\end{definition*}

The maximal singularities method allows to check $G$-birational (super)rigidity
using only internal geometry of the original variety, without considering all other $G$-Mfs's. We need the following technical definition which has become common nowadays.

\begin{definition*}
Let $X$ be a normal variety, let $\MMM$ be a linear system of Weil divisors on $X$ without fixed components, and let $\lambda$ be a rational number. 
We say that the pair $(X,\, \lambda \MMM)$ is \emph{canonical} if 
some multiple $m(K_X+\lambda M)$ is Cartier, where $M\in \MMM$, and
for any birational morphism
$f: Y\to X$ one can write
\[
m(K_Y+\lambda \MMM_Y) =f^*m(K_X+\lambda \MMM)+\sum a_i E_i, 
\]
where $\MMM_Y$ is the birational transform of $\MMM$, $E_i$ are prime exceptional divisors, and $a_i\ge 0$ for all $i$.
\end{definition*}

In the surface case the canonical property is very easy to check:
a pair $(X,\, \lambda \MMM)$ is canonical if and only if 
\[
\mult_P(\MMM)\le 1/\lambda
\]
for any point $P\in X$.

Now, suppose that a $\GQ$-Fano variety $X$ is not $G$-birationally superrigid. 
Then the Noether-Fano inequality \cite[Theorem 4.2]{Corti95:Sark} implies
the existence of a $G$-invariant linear system $\MMM$ on $X$ without fixed components such that
the pair $(X, \lambda \MMM)$ is not canonical, where $\lambda\in \QQ$ is taken so that $K_X + \lambda \MMM$ is numerically trivial. Moreover, any $\MMM$ as above defines a birational $G$-map $X \dashrightarrow X^{\sharp}$ to
the total space of a $G$-Mfs\ $X^{\sharp}/Z^{\sharp}$.
To show existence or non-existence of such $\MMM$ one needs to 
analyze the geometry of the variety $X$ carefully.

\begin{example*}
Let $X$ be a del Pezzo surface of degree $1$. 
Assume that $X$ is a $G$-del Pezzo with respect to some group $G\subset \Aut(X)$.
This means that $G$ acts on $X$ so that $\rk \Pic(X)^G=1$. For example,
this holds for any subgroup $G\subset \Aut(X)$ containing the Bertini involution. Let $\MMM$ be a $G$-invariant linear subsystem without fixed components. 
Since $\Pic(X)^G=\ZZ\cdot K_X$, we have $\MMM\subset |-nK_X|$ for some $n>0$.
Suppose that the pair $(X, \frac1n \MMM)$ is not canonical. Then $\mult_P(\MMM)>n$. Since $\MMM$ has no fixed components,
\[
n^2=(-nK_X)^2= \MMM^2\ge (\mult_P(\MMM))^2>n^2.
\]
The contradiction shows that $X$ is $G$-birationally superrigid.
\end{example*}

Similar arguments show that any $G$-del Pezzo surface $X$ of degree $\le 3$ is 
$G$-birationally rigid. Moreover, it is $G$-birationally superrigid if and only if $G$ has no orbits of length $\le K_X^2-2$ on $X$. In particular, $\PSL_2(\bF_7)$-del Pezzo surface from Theorem~\ref{thm:simple2} is $G$-birationally superrigid.
\begin{example*}
All the $\GQ$-Fano threefolds from Theorem~\ref{thm:simple3} are $G$-birationally superrigid
\cite{Cheltsov-Shramov:5, Cheltsov-Shramov:P3, BCDP:QS}. In particular, different embeddings of $\PSp_4(\bF_3)$ and $\PSL_2(\bF_{11})$ are not conjugate in $\Cr_3(\Bbbk)$. 
\end{example*}

There is  another relevant and very important notion called  $G$-solidity \cite{CDK:toric}.
For Fano varieties without group action
this notion has been  introduced earlier by Shokurov \cite{Shokurov1988p} (who called solid Fano varieties primitive) and by Abban and Okada \cite{Ahmadinezhad-Okada}. 

\begin{definition*}[{\cite{CDK:toric}}]
A $G$-Fano variety $X$ is \emph{$G$-solid} if $X$ is not $G$-birational to a $G$-Mfs with positive dimensional base.
\end{definition*}

For example a $G$-del Pezzo surface $X$ of degree $4$ is 
$G$-solid if and only if $G$ has no fixed points on $X$ \cite[\S~8]{Dolgachev-Iskovskikh}. 

A part of the maximal singularities method is so-called Sarkisov program \cite{Corti95:Sark, HaconMcKernan:Sark}. It allows us to decompose any birational map between Mfs's into a composition of elementary ones.
Refer to \cite{Iskovskikh:Factorization-e} for an explicit description of this program in dimension two and to \cite{Cheltsov:2ineq} for examples and applications.

\section{Application: essential dimension}
The notion of the essential dimension of a finite group $G$, denoted by $\ed(G)$, was introduced by Buhler and Reichstein \cite{Buhler-Reichstein-1997}. Informally, $\ed(G)$ is the minimal number of algebraic parameters needed 
to describe a 
faithful representation. More precisely, given a faithful linear representation $V$ of $G$ viewed as a $G$-variety,
the \emph{essential dimension} $\ed(G, V)$ 
is the minimal value of $\dim(X)$, where $X$ is taken from 
the set of all $G$-varieties 
admitting dominant rational 
$G$-equivariant map $V 
\dashrightarrow X$.
It can be shown that 
$\ed(G,V)$ does not depend on $V$, so we can omit $V$ in the notation. It is easy to see that $\ed(G) = 1$ if and only if $G$ is cyclic or dihedral of order $2n$ where $n$ is odd.
Finite groups of essential dimension $\le 2$ have been classified \cite{Duncan2013}. 

The essential dimension of symmetric groups $\Sym_n$ is important because it is equal to the minimal number of parameters needed to describe general polynomial
of degree $n$ modulo Tschirnhaus transformations \cite{Buhler-Reichstein-1997}.
The values of $\ed(\Sym_n)$, as well as, of $\ed(\Alt_n)$ are known for $n\le 7$ and bounds exist
for any $n$:
\begin{theorem}[\cite{Buhler-Reichstein-1997}, \cite{Duncan2010}]
If $n\ge 6$, then 
\begin{eqnarray*}
n - 3\ge &\ed (\Sym_n) &\ge \lfloor n/2 \rfloor,
\\
\ed (\Sym_n) \ge& \ed (\Alt_n) &\ge
\begin{cases}
\frac n2 &\text{if $n$ is even,}
\\ 
2 \lfloor \frac {n+2}4\rfloor&\text{if $n$ is odd.}
\end{cases}
\end{eqnarray*}
\end{theorem}

In many cases the computations of $\ed(G)$ use the machinery of $G$-varieties. 
As an example, following Serre \cite{Serre-2008-2009} we show that $\ed(\Alt_6)=3$. Let $V$ be the standard six-dimensional permutation representation of $\Alt_6$. There exists an equivariant open
embedding $V\subset (\PP^1)^6$. On the other hand, the group $\PSL_2(\Bbbk)$ also acts on $(\PP^1)^6$ so that the two actions commute.
Hence we have a dominant rational $\Alt_6$-map 
\[
V\hookrightarrow (\PP^1)^6 \longrightarrow (\PP^1)^6/\PSL_2(\Bbbk),
\]
where $(\PP^1)^6/\PSL_2(\Bbbk)$ is a birational quotient.
Since $(\PP^1)^6/\PSL_2(\Bbbk)=3$, we have $\ed(\Alt_6)\le 3$.
Thus it is sufficient to show that $\ed(\Alt_6)$ is not equal to $2$. 
If so, there exists a dominant rational $G$-map $V 
\dashrightarrow X$ to a surface which must be rational.
According to Theorem~\ref{thm:simple2} we may assume that $X=\PP^2$.
But in this case a Sylow $3$-subgroup $S\subset \Alt_6$ is abelian and acts without fixed points on $\PP^2$. On the other hand, $S$ has a fixed point on $V$ and the same should be true for 
the image of any rational $S$-map to a projective variety \cite{Kollar-Szabo-2000}. 
Therefore $\ed(\Alt_6)=3$ as claimed.

Using similar arguments and the classification of embeddings of $\Alt_7$ to groups of birational transformations of rationally connected threefolds (Theorem~\ref{thm:simple3}) A. Duncan proved that
$\ed(\Alt_7)=\ed(\Sym_7)=4$ \cite{Duncan2010}.

Denote by $\mathrm{rdim} (G)$ (resp. $\mathrm{cdim} (G)$) the minimal dimension of faithful 
representations of $G$ (resp. the smallest $n$ such that $G$ is embeddable to $\Cr_n(\Bbbk)$).
It is immediately follows from the definition that
\[
\ed(G) \le \mathrm{rdim} (G).
\]
If $G$ is a $p$-group, then the equality holds $\ed(G) = \mathrm{rdim} (G)$ \cite{Karpenko-Merkurjev-2008}. In general, this equality fails but 
there is a bound in terms of Jordan constants:

\begin{theorem}[\cite{Reichstein:J}]
$\mathrm{rdim} (G)\le \ed(G)\cdot \jj(\ed(G))$,
where $\jj(n)$ is the Jordan constant.
\end{theorem}

I. Dolgachev conjectured that $\ed(G) \ge \mathrm{cdim} (G)$ (see \cite{Duncan-Reichstein}). It would be interesting to test this conjecture for the group $G=\PSL_2(\bF_{11})$. In fact, we have
\[ 
3\le \ed(\PSL_2(\bF_{11}))\le 4
\]
by 
Theorem~\ref{thm:simple2} and because the group $\PSL_2(\bF_{11})$ is simple and has a faithful $5$-dimensional representation. 
Assuming Dolgachev's conjecture, by Theorem~\ref{thm:simple3} we would have $\ed(\PSL_2(\bF_{11}))=4$.
But this is unknown. See \cite{Duncan-Reichstein} for interesting discussions.
The computation of the essential dimension of $\PSL_2(\bF_{11})$ should complete Beauville's classification of finite simple groups of essential dimension $\le 3$ \cite{Beauville:ed-simple}.  

\begin{acknowledgments}
The author would like to thank Alexander Duncan, Constantin Shramov, Yuri Tschinkel, Mikhail Zaidenberg, and the referee 
for helpful comments on the original version of this paper.
\par\noindent
The study has been funded within the framework of the HSE University Basic Research Program.
\end{acknowledgments}

\def\cprime{$'$}


\begin{thebibliography}{100}

\bibitem{Abramovich-Wang}
{\sc Abramovich, D., and Wang, J.}
\newblock Equivariant resolution of singularities in characteristic {$0$}.
\newblock {\em Math. Res. Lett. 4}, 2-3 (1997), 427--433.

\bibitem{Ahmadinezhad-Okada}
{\sc Ahmadinezhad, H., and Okada, T.}
\newblock Birationally rigid {P}faffian {F}ano $3$-folds.
\newblock {\em Algebr. Geom. 5}, 2 (2018), 160--199.

\bibitem{Alberich-Carraminana}
{\sc Alberich-Carrami\~{n}ana, M.}
\newblock {\em Geometry of the plane {C}remona maps}, vol.~1769 of {\em Lecture
  Notes in Mathematics}.
\newblock Springer-Verlag, Berlin, 2002.

\bibitem{Avilov:cubic:Eu}
{\sc Avilov, A.}
\newblock Automorphisms of singular three-dimensional cubic hypersurfaces.
\newblock {\em Eur. J. Math. 4}, 3 (2018), 761--777.

\bibitem{Avilov:cb}
{\sc Avilov, A.~A.}
\newblock Existence of standard models of conic fibrations over
  non-algebraically-closed fields.
\newblock {\em Sb. Math. 205}, 12 (2014), 1683--1695.

\bibitem{Avilov-V4e}
{\sc Avilov, A.~A.}
\newblock Automorphisms of threefolds that can be represented as an
  intersection of two quadrics.
\newblock {\em Sb. Math. 207}, 3 (2016), 315--330.

\bibitem{BandmanZarhin:17}
{\sc Bandman, T., and Zarhin, Y.~G.}
\newblock Jordan groups, conic bundles and abelian varieties.
\newblock {\em Algebr. Geom. 4}, 2 (2017), 229--246.

\bibitem{Bayle-Beauville-2000}
{\sc Bayle, L., and Beauville, A.}
\newblock {Birational involutions of {${\bf P}\sp 2$}}.
\newblock {\em Asian J. Math. 4}, 1 (2000), 11--17.
\newblock Kodaira's issue.

\bibitem{Beauville2007}
{\sc Beauville, A.}
\newblock {{$p$}-elementary subgroups of the {C}remona group}.
\newblock {\em J. Algebra 314}, 2 (2007), 553--564.

\bibitem{Beauville:ed-simple}
{\sc Beauville, A.}
\newblock Finite simple groups of small essential dimension.
\newblock In {\em Trends in contemporary mathematics.. Selected talks based on
  the presentations at the INdAM day, June 18, 2014}. Cham: Springer, 2014,
  pp.~221--228.

\bibitem{BCTSSD85}
{\sc Beauville, A., Colliot-Th{\'e}l{\`e}ne, J.-L., Sansuc, J.-J., and
  Swinnerton-Dyer, P.}
\newblock {Vari{\'e}t{\'e}s stablement rationnelles non rationnelles}.
\newblock {\em Ann. of Math. (2) 121}, 2 (1985), 283--318.

\bibitem{Bertini1877}
{\sc Bertini, E.}
\newblock Ricerche sulle trasformazioni univoche involutorie nel piano.
\newblock {\em Annali di Mat. Pura Appl. 8\/} (1877), 254--287.

\bibitem{Birkar-BAB}
{\sc Birkar, C.}
\newblock Singularities of linear systems and boundedness of {F}ano varieties.
\newblock {\em Ann. of Math. (2) 193}, 2 (2021), 347--405.

\bibitem{BCHM}
{\sc Birkar, C., Cascini, P., Hacon, C.~D., and McKernan, J.}
\newblock Existence of minimal models for varieties of log general type.
\newblock {\em J. Amer. Math. Soc. 23}, 2 (2010), 405--468.

\bibitem{Blanc:Lin}
{\sc Blanc, J.}
\newblock Linearisation of finite abelian subgroups of the {C}remona group of
  the plane.
\newblock {\em Groups Geom. Dyn. 3}, 2 (2009), 215--266.

\bibitem{Blanc:cycl}
{\sc {Blanc}, J.}
\newblock Elements and cyclic subgroups of finite order of the {C}remona group.
\newblock {\em Comment. Math. Helv. 86}, 2 (2011), 469--497.

\bibitem{BCDP:QS}
{\sc Blanc, J., Cheltsov, I., Duncan, A., and Prokhorov, Y.}
\newblock Finite quasisimple groups acting on rationally connected threefolds.
\newblock {\em ArXiv e-print 1809.09226\/} (2018).

\bibitem{BlancLZ:quotients}
{\sc {Blanc}, J., {Lamy}, S., and {Zimmermann}, S.}
\newblock Quotients of higher-dimensional cremona groups.
\newblock {\em Acta Math. 226}, 2 (2021), 211--318.

\bibitem{Bogomolov-Prokhorov}
{\sc Bogomolov, F., and Prokhorov, Y.}
\newblock {On stable conjugacy of finite subgroups of the plane {C}remona
  group, {I}}.
\newblock {\em Cent. European J. Math. 11}, 12 (2013), 2099--2105.

\bibitem{Braun-Filipazzi-Moraga-Svaldi}
{\sc Braun, L., Filipazzi, S., Moraga, J., and Svaldi, R.}
\newblock The {J}ordan property for local fundamental groups.
\newblock {\em Arxiv e-print 2006.01253\/} (2020).

\bibitem{Buhler-Reichstein-1997}
{\sc Buhler, J., and Reichstein, Z.}
\newblock On the essential dimension of a finite group.
\newblock {\em Compos. Math. 106}, 2 (1997), 159--179.

\bibitem{Cantat2013}
{\sc Cantat, S.}
\newblock The {C}remona group in two variables.
\newblock In {\em European {C}ongress of {M}athematics}. Eur. Math. Soc.,
  Z\"{u}rich, 2013, pp.~211--225.

\bibitem{Cantat2018}
{\sc Cantat, S.}
\newblock The {C}remona group.
\newblock In {\em Algebraic geometry: {S}alt {L}ake {C}ity 2015}, vol.~97 of
  {\em Proc. Sympos. Pure Math.} Amer. Math. Soc., Providence, RI, 2018,
  pp.~101--142.

\bibitem{CantatLamy}
{\sc Cantat, S., and Lamy, S.}
\newblock Normal subgroups in the {C}remona group.
\newblock {\em Acta Math. 210}, 1 (2013), 31--94.
\newblock With an appendix by Yves de Cornulier.

\bibitem{CDK:toric}
{\sc Cheltsov, I., Dubouloz, A., and Kishimoto, T.}
\newblock Toric {$G$}-solid {F}ano threefolds.
\newblock {\em Arxiv e-print 2007.14197\/} (2020).

\bibitem{Cheltsov-Przyjalkowski-Shramov-Barth-sextic}
{\sc Cheltsov, I., Przyjalkowski, V., and Shramov, C.}
\newblock Burkhardt quartic, {B}arth sextic, and the icosahedron.
\newblock {\em Int. Math. Res. Not. IMRN}, 12 (2019), 3683--3703.

\bibitem{Cheltsov-Shramov-2012}
{\sc Cheltsov, I., and Shramov, C.}
\newblock Three embeddings of the {K}lein simple group into the {C}remona group
  of rank three.
\newblock {\em Transform. Groups 17}, 2 (2012), 303--350.

\bibitem{Cheltsov-Shramov:5}
{\sc Cheltsov, I., and Shramov, C.}
\newblock Five embeddings of one simple group.
\newblock {\em Trans. Amer. Math. Soc. 366}, 3 (2014).

\bibitem{CheltsovShramov:book}
{\sc Cheltsov, I., and Shramov, C.}
\newblock {\em Cremona groups and the icosahedron}.
\newblock Boca Raton, FL: CRC Press, 2016.

\bibitem{Cheltsov-Shramov:P3}
{\sc Cheltsov, I., and Shramov, C.}
\newblock Finite collineation groups and birational rigidity.
\newblock {\em Sel. Math., New Ser. 25}, 5 (2019), 68.
\newblock Id/No 71.

\bibitem{Cheltsov:2ineq}
{\sc Cheltsov, I.~A.}
\newblock Two local inequalities.
\newblock {\em Izv. Math. 78}, 2 (2014), 375--426.

\bibitem{Conrad:Lang-Neron}
{\sc Conrad, B.}
\newblock Chow's {$K/k$}-image and {$K/k$}-trace, and the {Lang-N\'eron}
  theorem.
\newblock {\em Enseign. Math. (2) 52}, 1-2 (2006), 37--108.

\bibitem{atlas}
{\sc Conway, J., Curtis, R., Norton, S., Parker, R., and Wilson, R.}
\newblock {\em {Atlas of finite groups. {M}aximal subgroups and ordinary
  characters for simple groups. {W}ith comput. assist. from {J}. {G}.
  {T}hackray}}.
\newblock Oxford: Clarendon Press. XXXIII, 252 p., 1985.

\bibitem{Corti95:Sark}
{\sc Corti, A.}
\newblock Factoring birational maps of threefolds after {S}arkisov.
\newblock {\em J. Algebraic Geom. 4}, 2 (1995), 223--254.

\bibitem{Corti:DP}
{\sc Corti, A.}
\newblock Del {P}ezzo surfaces over {D}edekind schemes.
\newblock {\em Ann. of Math. (2) 144}, 3 (1996), 641--683.

\bibitem{Fernex:Crem}
{\sc de~Fernex, T.}
\newblock On planar {C}remona maps of prime order.
\newblock {\em Nagoya Math. J. 174\/} (2004), 1--28.

\bibitem{Debarre-Kuznetsov:GM}
{\sc Debarre, O., and Kuznetsov, A.}
\newblock Gushel-{M}ukai varieties: classification and birationalities.
\newblock {\em Algebr. Geom. 5}, 1 (2018), 15--76.

\bibitem{Deserti:book-s}
{\sc D{\'e}serti, J.}
\newblock {\em Some properties of the {C}remona group}, vol.~21 of {\em Ensaios
  Matem{\'a}ticos [Mathematical Surveys]}.
\newblock Sociedade Brasileira de Matem{\'a}tica, Rio de Janeiro, 2012.

\bibitem{Dolgachev:fields-inv}
{\sc Dolgachev, I.~V.}
\newblock Rationality of fields of invariants.
\newblock In {\em Algebraic geometry, {B}owdoin, 1985 ({B}runswick, {M}aine,
  1985)}, vol.~46 of {\em Proc. Sympos. Pure Math.} Amer. Math. Soc.,
  Providence, RI, 1987, pp.~3--16.

\bibitem{Dolgachev-Iskovskikh}
{\sc Dolgachev, I.~V., and Iskovskikh, V.~A.}
\newblock Finite subgroups of the plane {C}remona group.
\newblock In {\em {Algebra, arithmetic, and geometry: in honor of {Y}u. {I}.
  {M}anin. {V}ol. {I}}}, vol.~269 of {\em {Progr. Math.}} Birkh{\"a}user Boston
  Inc., Boston, MA, 2009, pp.~443--548.

\bibitem{Dolgachev-Iskovskikh:p}
{\sc Dolgachev, I.~V., and Iskovskikh, V.~A.}
\newblock {On elements of prime order in the plane {C}remona group over a
  perfect field}.
\newblock {\em Int. Math. Res. Not. IMRN}, 18 (2009), 3467--3485.

\bibitem{Duncan2010}
{\sc Duncan, A.}
\newblock Essential dimensions of {$A_7$} and {$S_7$}.
\newblock {\em Math. Res. Lett. 17}, 2 (2010), 263--266.

\bibitem{Duncan2013}
{\sc Duncan, A.}
\newblock Finite groups of essential dimension 2.
\newblock {\em Comment. Math. Helv. 88}, 3 (2013), 555--585.

\bibitem{Duncan-Reichstein}
{\sc Duncan, A., and Reichstein, Z.}
\newblock Versality of algebraic group actions and rational points on twisted
  varieties.
\newblock {\em J. Algebraic Geom. 24}, 3 (2015), 499--530.
\newblock With an appendix containing a letter from J.-P. Serre.

\bibitem{HaconMcKernan:Sark}
{\sc Hacon, C.~D., and McKernan, J.}
\newblock The {S}arkisov program.
\newblock {\em J. Algebraic Geom. 22}, 2 (2013), 389--405.

\bibitem{Hassett-Kresch-Tschinkel:symbols}
{\sc Hassett, B., Kresch, A., and Tschinkel, Y.}
\newblock Symbols and equivariant birational geometry in small dimensions.
\newblock {\em Arxiv e-rint\/} (2020).

\bibitem{Haution:19}
{\sc Haution, O.}
\newblock Fixed point theorems involving numerical invariants.
\newblock {\em Compos. Math. 155}, 2 (2019), 260–--288.

\bibitem{Hudson1927}
{\sc Hudson, H.~P.}
\newblock {\em {Cremona transformations in plane and space}}.
\newblock Cambridge University Press, 1927.

\bibitem{Isk-Manin}
{\sc Iskovskikh, V., and Manin, Y.}
\newblock Three-dimensional quartics and counterexamples to the {L{\"u}}roth
  problem.
\newblock {\em Math. USSR-Sb. 15}, 1 (1971), 141--166.

\bibitem{Iskovskikh:Factorization-e}
{\sc Iskovskikh, V.~A.}
\newblock Factorization of birational mappings of rational surfaces from the
  point of view of {M}ori theory.
\newblock {\em Russian Math. Surveys 51}, 4 (1996), 585--652.

\bibitem{Isk:Two1e}
{\sc Iskovskikh, V.~A.}
\newblock Two nonconjugate embeddings of the group {\(s_3\times \mathbb z_2\)
  into the cremona group}.
\newblock {\em Proc. Steklov Inst. Math. 241\/} (2003), 93--97.

\bibitem{IP99}
{\sc Iskovskikh, V.~A., and Prokhorov, Y.}
\newblock {\em Fano varieties. {A}lgebraic geometry {V}}, vol.~47 of {\em
  {Encyclopaedia Math. Sci.}}
\newblock Springer, Berlin, 1999.

\bibitem{Jordan:1878}
{\sc Jordan, C.}
\newblock M\'emoire sur les \'equations diff\'erentielles lin\'eaires \`a
  int\'egrale alg\'ebrique.
\newblock {\em {J. Reine Angew. Math.} 84\/} (1878), 89--215.

\bibitem{Karpenko-Merkurjev-2008}
{\sc Karpenko, N., and Merkurjev, A.}
\newblock {Essential dimension of finite $p$-groups}.
\newblock {\em Invent. math. 172\/} (2008), 491--508.

\bibitem{Kawamata:bF}
{\sc Kawamata, Y.}
\newblock Boundedness of {$\mathbf{Q}$}-{F}ano threefolds.
\newblock In {\em Proceedings of the {I}nternational {C}onference on {A}lgebra,
  {P}art 3 ({N}ovosibirsk, 1989)\/} (1992), vol.~131 of {\em Contemp. Math.},
  Amer. Math. Soc., Providence, RI, pp.~439--445.

\bibitem{Kollar-Miyaoka-Mori-1992a}
{\sc Koll{\'a}r, J., Miyaoka, Y., and Mori, S.}
\newblock {Rationally connected varieties}.
\newblock {\em J. Algebraic Geom. 1}, 3 (1992), 429--448.

\bibitem{KMMT-2000}
{\sc Koll{\'a}r, J., Miyaoka, Y., Mori, S., and Takagi, H.}
\newblock Boundedness of canonical {$\mathbf{Q}$}-{F}ano 3-folds.
\newblock {\em Proc. Japan Acad. Ser. A Math. Sci. 76}, 5 (2000), 73--77.

\bibitem{Kollar-Szabo-2000}
{\sc Koll{\'a}r, J., and Szab{\'o}, E.}
\newblock Fixed points of group actions and rational maps.
\newblock {\em Canad. J. Math. 52}, 5 (2000), 1054--1056.
\newblock {Appendix to ``Essential dimensions of algebraic groups and a
  resolution theorem for {$G$}-varieties'' by Z. Reichstein and B. Youssin}.

\bibitem{kontsevich-Pestun-Tschinkel}
{\sc Kontsevich, M., Pestun, V., and Tschinkel, Y.}
\newblock Equivariant birational geometry and modular symbols.
\newblock {\em Arxiv e-print 1902.09894\/} (2019).

\bibitem{Kresch-Tschinkel:Burnside-vol}
{\sc Kresch, A., and Tschinkel, Y.}
\newblock Equivariant birational types and {B}urnside volume.
\newblock {\em Arxiv e-print 2007.12538\/} (2020).
\newblock to appear in Ann. Scuola Norm. Sup. Pisa.

\bibitem{Kresch-Tschinkel:Burnside-repr}
{\sc Kresch, A., and Tschinkel, Y.}
\newblock Equivariant burnside groups and representation theory.
\newblock {\em Arxiv e-print 2108.00518\/} (2021).

\bibitem{Krylov2020FamiliesOE}
{\sc Krylov, I.}
\newblock Families of embeddings of the alternating group of rank $5$ into the
  {C}remona group.
\newblock {\em Arxiv e-print 2005.07354\/} (2020).

\bibitem{KPS:Hilb}
{\sc Kuznetsov, A., Prokhorov, Y., and Shramov, C.}
\newblock Hilbert schemes of lines and conics and automorphism groups of {F}ano
  threefolds.
\newblock {\em Japanese J. Math. 13}, 1 (2018), 109--185.

\bibitem{Kuznetsova:Finite3}
{\sc Kuznetsova, A.}
\newblock Finite $3$-subgroups in {C}remona group of rank $3$.
\newblock {\em Math. Notes 108\/} (2020), 697--715.

\bibitem{Lemire-Popov-Reichstein}
{\sc Lemire, N., Popov, V.~L., and Reichstein, Z.}
\newblock Cayley groups.
\newblock {\em J. Amer. Math. Soc. 19}, 4 (2006), 921--967.

\bibitem{Loginov:dP1}
{\sc Loginov, K.}
\newblock Standard models of degree $1$ del {P}ezzo fibrations.
\newblock {\em Mosc. Math. J. 18}, 4 (2018), 721--737.

\bibitem{Loginov:3Cr}
{\sc Loginov, K.}
\newblock A note on $3$-subgroups in the space {C}remona group.
\newblock {\em ArXiV e-print 2102.04522\/} (2021).

\bibitem{Manin:raf-surf2e}
{\sc Manin, Y.~I.}
\newblock {Rational surfaces over perfect fields. {II}}.
\newblock {\em Mat. Sb. (N.S.) 72 (114)\/} (1967), 161--192.

\bibitem{Minkowski:87}
{\sc Minkowski, H.}
\newblock Zur {T}heorie der positiven quadratischen {F}ormen.
\newblock {\em J. Reine Angew. Math. 101\/} (1887), 196--202.

\bibitem{MiyaokaMori}
{\sc Miyaoka, Y., and Mori, S.}
\newblock {A numerical criterion for uniruledness}.
\newblock {\em Ann. of Math. (2) 124}, 1 (1986), 65--69.

\bibitem{Moraga:toroidalization}
{\sc Moraga, J.}
\newblock On a toroidalization for klt singularities.
\newblock {\em Arxiv e-print 2106.15019\/} (2021).

\bibitem{Mori:term-sing}
{\sc Mori, S.}
\newblock On {$3$}-dimensional terminal singularities.
\newblock {\em Nagoya Math. J. 98\/} (1985), 43--66.

\bibitem{Mori:flip}
{\sc Mori, S.}
\newblock Flip theorem and the existence of minimal models for {$3$}-folds.
\newblock {\em J. Amer. Math. Soc. 1}, 1 (1988), 117--253.

\bibitem{Namikawa:Fano}
{\sc Namikawa, Y.}
\newblock Smoothing {F}ano {$3$}-folds.
\newblock {\em J. Algebraic Geom. 6}, 2 (1997), 307--324.

\bibitem{Popov:Russel}
{\sc Popov, V.~L.}
\newblock On the {M}akar-{L}imanov, {D}erksen invariants, and finite
  automorphism groups of algebraic varieties.
\newblock In {\em Peter Russell's Festschrift, Proceedings of the conference on
  Affine Algebraic Geometry held in Professor Russell's honour, 1--5 June 2009,
  McGill Univ., Montreal\/} (2011), vol.~54 of {\em {Centre de Recherches
  Math{\'e}matiques CRM Proc. and Lect. Notes}}, pp.~289--311.

\bibitem{Popov-tori}
{\sc Popov, V.~L.}
\newblock Tori in the {C}remona groups.
\newblock {\em Izv. Math. 77}, 4 (2013), 742--771.

\bibitem{P:p-groups}
{\sc Prokhorov, Y.}
\newblock {$p$}-elementary subgroups of the {C}remona group of rank $3$.
\newblock In {\em Classification of algebraic varieties}, {EMS Ser. Congr.
  Rep.} Eur. Math. Soc., Z{\"u}rich, 2011, pp.~327--338.

\bibitem{P:JAG:simple}
{\sc Prokhorov, Y.}
\newblock Simple finite subgroups of the {C}remona group of rank $3$.
\newblock {\em J. Algebraic Geom. 21}, 3 (2012), 563--600.

\bibitem{P:GFano-all}
{\sc Prokhorov, Y.}
\newblock G-{F}ano threefolds, {I}, {II}.
\newblock {\em Adv. Geom. 13}, 3 (2013), 389--418, 419--434.

\bibitem{P:invol}
{\sc Prokhorov, Y.}
\newblock {On birational involutions of {$\mathbf P^3$}}.
\newblock {\em Izvestiya: Math. 77}, 3 (2013), 627--648.

\bibitem{Prokhorov-2-elementary}
{\sc Prokhorov, Y.}
\newblock {2-elementary subgroups of the space {C}remona group}.
\newblock In {\em {Automorphisms in birational and affine geometry}}, vol.~79
  of {\em {Springer Proc. Math. Stat.}} Springer, Cham, 2014, pp.~215--229.

\bibitem{Prokhorov-v22}
{\sc Prokhorov, Y.}
\newblock Singular {F}ano threefolds of genus $12$.
\newblock {\em Sbornik: Math. 207}, 7 (2016), 983--1009.

\bibitem{P:G-MMP}
{\sc Prokhorov, Y.}
\newblock Equivariant minimal model program.
\newblock {\em Russian Math. Surv. 76}, 3 (2021), 461--542.

\bibitem{Prokhorov-Shramov-J}
{\sc Prokhorov, Y., and Shramov, C.}
\newblock Jordan property for groups of birational selfmaps.
\newblock {\em Compositio Math. 150}, 12 (2014), 2054--2072.

\bibitem{ProkhorovShramov-RC}
{\sc Prokhorov, Y., and Shramov, C.}
\newblock Jordan property for {C}remona groups.
\newblock {\em Amer. J. Math. 138}, 2 (2016), 403--418.

\bibitem{Prokhorov-Shramov-J-const}
{\sc Prokhorov, Y., and Shramov, C.}
\newblock Jordan constant for {C}remona group of rank {$3$}.
\newblock {\em Moscow Math. J. 17}, 3 (2017), 457--509.

\bibitem{P-Shramov:3folds}
{\sc Prokhorov, Y., and Shramov, C.}
\newblock Finite groups of birational selfmaps of threefolds.
\newblock {\em Math. Res. Lett. 25}, 3 (2018), 957--972.

\bibitem{Prokhorov-Shramov-p-groups}
{\sc Prokhorov, Y., and Shramov, C.}
\newblock {$p$}-subgroups in the space {C}remona group.
\newblock {\em Math. Nachrichten 291}, 8--9 (2018), 1374--1389.

\bibitem{Pukhlikov:book}
{\sc Pukhlikov, A.}
\newblock {\em Birationally rigid varieties}, vol.~190 of {\em Mathematical
  Surveys and Monographs}.
\newblock American Mathematical Society, Providence, RI, 2013.

\bibitem{Reichstein:J}
{\sc Reichstein, Z.}
\newblock The {J}ordan property of {C}remona groups and essential dimension.
\newblock {\em Arch. Math. 111}, 5 (2018), 449--455.

\bibitem{Reichstein-Youssin:bir}
{\sc Reichstein, Z., and Youssin, B.}
\newblock A birational invariant for algebraic group actions.
\newblock {\em Pacific J. Math. 204}, 1 (2002), 223--246.

\bibitem{Reid:YPG}
{\sc Reid, M.}
\newblock Young person's guide to canonical singularities.
\newblock In {\em Algebraic geometry, {B}owdoin, 1985 ({B}runswick, {M}aine,
  1985)}, vol.~46 of {\em Proc. Sympos. Pure Math.} Amer. Math. Soc.,
  Providence, RI, 1987, pp.~345--414.

\bibitem{Serre:MMJ}
{\sc Serre, J.-P.}
\newblock {A {M}inkowski-style bound for the orders of the finite subgroups of
  the {C}remona group of rank 2 over an arbitrary field}.
\newblock {\em Mosc. Math. J. 9}, 1 (2009), 193--208.

\bibitem{Serre-2008-2009}
{\sc {Serre}, J.-P.}
\newblock {Le groupe de {C}remona et ses sous-groupes finis}.
\newblock In {\em {S{\'e}minaire Bourbaki. Volume 2008/2009. Expos{\'e}s
  997--1011}}. Paris: Soci{\'e}t{\'e} Math{\'e}matique de France (SMF), 2010,
  pp.~75--100, ex.

\bibitem{Serre:problems}
{\sc Serre, J.-P.}
\newblock Problems for the {E}dinburgh workshop on {C}remona groups, March
  2010.

\bibitem{Shinder}
{\sc Shinder, E.}
\newblock The {B}ogomolov-{P}rokhorov invariant of surfaces as equivariant
  cohomology.
\newblock {\em Bull. Korean Math. Soc. 54}, 5 (2017), 1725--1741.

\bibitem{Shokurov1988p}
{\sc Shokurov, V.~V.}
\newblock {Problems about {F}ano varieties}.
\newblock In {\em {Birational Geometry of Algebraic Varieties, Open
  Problems}\/} (Katata, 1988), pp.~30--32.

\bibitem{Shokurov:PL:e}
{\sc Shokurov, V.~V.}
\newblock Prelimiting flips.
\newblock {\em Proc. Steklov Inst. Math. 240\/} (2003), 75--213.

\bibitem{Shramov:aut-2cubics}
{\sc Shramov, C.}
\newblock Automorphisms of cubic surfaces without points.
\newblock {\em Int. J. Math. 31}, 11 (2020), 17.
\newblock Id/No 2050083.

\bibitem{Shramov:SB-bir}
{\sc Shramov, C.}
\newblock Birational automorphisms of {S}everi-{B}rauer surfaces.
\newblock {\em Sb. Math. 211}, 3 (2020), 466--480.

\bibitem{Shramov:SBf}
{\sc Shramov, C.}
\newblock Finite groups acting on {S}everi-{B}rauer surfaces.
\newblock {\em Eur. J. Math. 7}, 2 (2021), 591--612.

\bibitem{Tsygankov:MSb:e}
{\sc Tsygankov, V.~I.}
\newblock Equations of {$G$}-minimal conic bundles.
\newblock {\em Sb. Math. 202}, 11-12 (2011), 1667--1721.
\newblock translated from Mat. Sb. 202 (2011), no. 11, 103--160.

\bibitem{Xu:p-groups}
{\sc Xu, J.}
\newblock A remark on the rank of finite $p$-groups of birational
  automorphisms.
\newblock {\em Comptes Rendus. Math\'ematique 358}, 7 (2020), 827--829.

\bibitem{Yasinsky-J-const}
{\sc Yasinsky, E.}
\newblock The {J}ordan constant for {C}remona group of rank {$2$}.
\newblock {\em Bull. Korean Math. Soc. 54}, 5 (2017), 1859--1871.

\bibitem{Zarhin10}
{\sc Zarhin, Y.~G.}
\newblock {Theta groups and products of abelian and rational varieties}.
\newblock {\em Proc. Edinburgh Math. Soc. 57}, 1 (2014), 299--304.

\end{thebibliography}

\end{document}